\documentclass[10pt,conference]{IEEEtran}
\usepackage[letterpaper, left=0.625in, right=0.625in, bottom=0.9in, top=0.7in]{geometry}
\IEEEoverridecommandlockouts
\usepackage[noadjust]{cite}
\usepackage{cite}
\usepackage{amsmath,amssymb,amsfonts,dsfont}
\usepackage{algorithm} 
\usepackage{algpseudocode} 
\usepackage{graphicx}
\usepackage{textcomp}
\usepackage{xcolor}
\usepackage{comment}
\usepackage{adjustbox}
\usepackage{float}
\usepackage[caption=false,labelformat=empty]{subfig}

\DeclareMathOperator*{\argmax}{arg\,max}
\DeclareMathOperator*{\argmin}{arg\,min}

\newcommand{\E}{\mathbb{E}}
\newcommand{\R}{\mathbb{R}}
\newcommand{\N}{\mathbb{N}}
\newcommand{\Z}{\mathbb{Z}}
\newcommand{\Char}{\mathds{1}}

\newcommand{\cS}{\mathcal{S}}

\newcommand{\cN}{\mathcal{N}}

\newcommand{\floor}[1]{\!\left\lfloor #1 \right\rfloor\!}

\newcommand{\pibf}{\boldsymbol{\pi}}

\newtheorem{theorem}{Theorem}[section]

\newtheorem{corollary}[theorem]{Corollary}
\newtheorem{proposition}[theorem]{Proposition}

\begin{document}

\title{Opportunistic Spectrum Access: Does Maximizing Throughput Minimize File Transfer Time?
\thanks{{\rule{5cm}{2pt}}

Jie Hu and Do Young Eun are with the Department of Electrical and Computer Engineering, and Vishwaraj Doshi is with the Operations Research Graduate Program, North Carolina State University, Raleigh, NC. Email: \{jhu29, vdoshi, dyeun\}@ncsu.edu. This work was supported in part by National Science Foundation under Grant CNS-1824518 and IIS-1910749.}
}
\author{\IEEEauthorblockN{Jie Hu}
\and
\IEEEauthorblockN{Vishwaraj Doshi}
\and
\IEEEauthorblockN{Do Young Eun}
}
\maketitle
\begin{abstract}
The \textit{Opportunistic Spectrum Access} (OSA) model has been developed for the secondary users (SUs) to exploit the stochastic dynamics of licensed channels for file transfer in an opportunistic manner. Common approaches to design channel sensing strategies for throughput-oriented applications tend to maximize the long-term throughput, with the hope that it provides reduced file transfer time as well. In this paper, we show that this is not correct in general, especially for small files. Unlike prior delay-related works that seldom consider the heterogeneous channel rate and bursty incoming packets, our work explicitly considers minimizing the file transfer time of a \emph{single} file consisting of multiple packets in a set of heterogeneous channels. We formulate a mathematical framework for the \emph{static} policy, and extend to \emph{dynamic} policy by mapping our file transfer problem to the stochastic shortest path problem. We analyze the performance of our proposed static optimal and dynamic optimal policies over the policy that maximizes long-term throughput. We then propose a heuristic policy that takes into account the performance-complexity tradeoff and an extension to online implementation with unknown channel parameters, and also present the regret bound for our online algorithm. We also present numerical simulations that reflect our analytical results.
\end{abstract}

\section{Introduction}\label{Introduction}

In recent years, rapid growth in the number of wireless devices, including mobile devices and Internet of Things (IoT) devices has led to an explosion in the demand for wireless service. This demand further exacerbates the scarcity of allocated spectrum, which is ironically known to be underutilized by licensed users \cite{osa_survey}. The \textit{Opportunistic spectrum access} (OSA) model has been proposed to reuse the licensed spectrum in an opportunistic way otherwise wasted by licensed users \cite{osa_survey}. Recently, the FCC has released a new guidance in $2020$, which would expand the ability of the unlicensed devices (especially IoT devices) to operate in the TV-broadcast bands \cite{federal_2020}. Besides, the related IEEE 802.22 family has been developed to enable spectrum sharing \cite{802.22_standard} to bring broadband access to rural areas.

In the OSA model, a \textit{secondary user} (SU) aims to opportunistically access the spectrum when it is not used by any other users, while also prioritizing the needs of the \textit{primary user} (PU). The SUs need to periodically sense the spectrum to avoid interfering with PUs. Interference reduces the quality of service in the OSA networks, and \emph{throughput} is one of the most commonly used performance metrics in the OSA literature. The PU's behavior in a channel can be modeled as a Markov process (thus correlated over time), for which partially observable Markov decision processes (POMDPs) are typically employed to formulate the spectrum sensing strategy in order to maximize the long-term throughput \cite{Zhao_2007}. These POMDPs do not possess known structured solutions in general and they are known to be PSPACE-complete even if all the channel statistics are known a priori \cite{liu2015online}. To achieve near maximum throughput, computationally efficient yet sub-optimal policies, such as myopic policy \cite{zhao2008opportunistic} and Whittle's index policy \cite{liu2010indexability}, have been proposed for the \emph{offline} OSA setting (known channel parameters) and later extended to the \emph{online} setting (unknown channel parameters) \cite{NEURIPS2019_2edfeadf}. Besides, single-channel online policies have also been developed in \cite{tekin2011online,dai2012efficient} to find the channel with maximum throughput in the steady state. Recently, by focusing more on heterogeneous channels with \emph{i.i.d} Bernoulli distribution for each, multi-armed bandit (MAB) techniques have been extensively studied to let SUs learn unknown channel parameters on the fly, and have been known for their \emph{lightweight} designs. MAB techniques that maximize the cumulative reward (or minimize the regret) can directly be applied to the online setting for maximizing the throughput, including the Bayesian approach \cite{Poor_2011}, upper confidence bounds \cite{Optimal_802.11}, thompson sampling \cite{Srikant_2018} and its improvement from efficient sampling \cite{Srikant_2019}, and coordination approach among multiple SUs \cite{avner2019multi}.

Nowadays, low latency has become one of the main goals for $5$G wireless networks \cite{Low_Latency_Towards_5G} and other time-sensitive applications with guaranteed delay constraints.
For example, packet delay in cognitive radio networks has been extensively studied using queuing theory to derive delay-efficient spectrum scheduling strategies. In this setting, a stream of packet arrivals modeled as a Poisson process with a constant rate is a common assumption in delay related works \cite{shiang2008queuing, wu2014learning, dimitriou2018stable,cao2017dynamic, huang2019dynamic}, and the goal is often to minimize the average packet delay in the steady state. In reality, however, so-called `bursty' cases are also commonly seen, where a finite number of packets comprising one file can be pushed into the SU's queue simultaneously and transmitted opportunistically by the SU, and the next `file' will not arrive if the current file transfer job has not been finished yet.\footnote{One packet transmission in the delay-related studies \cite{shiang2008queuing,wu2014learning,cao2017dynamic,dimitriou2018stable,huang2019dynamic} is equivalent to a single file transfer job in the `bursty' case because a single file in the SU's queue can be seen as a `large' packet.} For instance, the IEEE 802.11p MAC protocol requires the safety message to be generated and sent by each vehicle in every $100$ ms interval \cite{8024608}, Poisson arrivals being unsuitable to model such situation. Same channel rate across all channels is another implicit assumption in many works \cite{wu2014learning,cao2017dynamic,dimitriou2018stable,huang2019dynamic}, but it doesn't reflect the realistic heterogeneous channel environment frequently assumed in the throughput-oriented studies in the OSA literature \cite{Poor_2011,Optimal_802.11,Srikant_2018,Srikant_2019,dai2012efficient}. Clearly, allowing the SU to switch across heterogeneous-rate channels during instances of PU’s interruption can further reduce the file transfer time, but to the best of our knowledge, this issue has not been fully explored.\footnote{Although \cite{shiang2008queuing} assumes heterogeneous channel rates, it requires the SU to stick to one channel to complete the packet transmission, even if the transmission may be interrupted by the PU multiple times.}

In this paper, we study the OSA model with the aim of minimizing the transfer time of a single file over heterogeneous channels. The common folklore assumes that the policy maximizing the long-term throughput would also lead to the minimum expected file transfer time. For example, Wald's equation implies that the file download time in an \emph{i.i.d} (over time) channel is equal to the file size divided by the average throughput of that channel, implicitly favoring the max-throughput channel for minimal download time. In this paper, we show that this is \emph{not} the case in general, even in the \emph{i.i.d} (over time) channel. We first present a theoretical analysis for \textit{static} policies where the SU only sticks to one channel throughout the entire file transfer. By casting the file transfer problem into a stochastic shortest path framework, the SU is free to switch between channels and we are able to obtain the dynamic \textit{optimal} policy. Our theoretical analysis shows that static and dynamic optimal policies reduce the transfer time compared to the baseline max-throughput policy, and this reduction is even more significant in delay-sensitive applications, where files are relatively small. We then propose a lightweight heuristic policy with good performance and extend to an online implementation with unknown channel parameters the SU needs to learn on the fly. We modify a MAB algorithm (model-based method) proposed in \cite{talebi2017stochastic} in the online setting and show its gap-dependent regret bound, instead of the model-free reinforcement learning used in delay-related works \cite{wu2014learning,cao2017dynamic} for sample efficiency purpose \cite{tu2019gap}. We also use simulation results to visualize that the max-throughput policy is not the best when it comes to achieving the minimum file transfer time.

The rest of the paper is organized as follows: In Section \ref{model description} we introduce the OSA model and characterize the file transfer problem and its policy under the OSA framework. In Section \ref{static analysis}, we show the expected transfer time for static policy and it's performance analysis. Then, we extend from the static policy to the dynamic policy in Section \ref{mdp formulation and dynamic policy}. The practical concerns are discussed in Section \ref{online learning}. Finally in Section \ref{simulation}, we evaluate different policies in the realistic numerical setting.

\section{Model description}\label{model description}
\subsection{The OSA Model}\label{problem setting}
Consider a set of $N$ heterogeneous channels $\cN \triangleq \{1,2,\cdots,N\}$ available for use and each channel $i\in\cN$ offers a stable rate of $r_i>0$ bits/s if successfully utilized \cite{dai2012efficient,liu2010indexability}. In our setting, a SU wishes to transfer a file of size $F$ bits using one of these $N$ channels via opportunistic spectrum access. The SU can only access one channel at any given time, and can maintain this access for a fixed duration of $\Delta$ seconds, after which it has to sense available channels again in order to access them. At this point, the SU can decide which channel to sense and access that channel for the next $\Delta$ second interval if the channel is \textit{available}. Or the SU has to wait for $\Delta$ seconds to sense again if that channel is \textit{unavailable}, thereby unable to transfer data for this duration. This pattern is known as the \textit{constant access time} model, and has been commonly adopted for the SU's behavior as a collision prevention mechanism in the OSA literature \cite{osa_survey,Poor_2011}.
The cycle repeats itself untill the SU transmits the entire file size $F$, then it immediately exits the channel in use. Note that the duration $\Delta$ seconds is not a randomly chosen number. For example, $\Delta$ is recommended as $100$ ms because the SU needs to vacate the current channel within $100$ ms once the PU shows up, as defined in IEEE 802.22 standard \cite{802.22_standard}. The SU can transmit up to $3.1$ Mb in each $\Delta$ seconds with highest channel rate $31$ Mbps in IEEE 802.22 standard and many small files (e.g, $5$ KB text-only email, $800$ KB GIF image) need just a few slots to complete.

We say a channel is unavailable (or \textit{busy}) if it is currently in use by the primary users (PUs) or other SUs, while it is available (or \textit{idle}) if it is not in use by any other users. The state of a channel (idle or busy) is assumed to be independent over all channels $i\in\cN$, and \emph{i.i.d.} over the time instants $\{0,\Delta,2\Delta,\cdots \}$ following Bernoulli distribution with parameter $p_i\in(0,1]$, in line with the widely used discrete-time channel model \cite{osa_survey}. Specifically, for each $i \in \cN$, $\{Y_i(k)\}_{k\in\N}$ is a Bernoulli process with $p_i=P[Y_i(k)=1]=1-P[Y_i(k)=0]$ for all $k\in\N$.
Then, we can define $X_i(t)$, the state of channel $i\in\cN$ at any time $t \in \R_+$, as a piecewise constant random process:
\begin{equation}\label{channel availability random variable}
    X_i(t) \triangleq Y_i\left(\floor{\frac{t}{\Delta}} \right),
\end{equation}
where $\floor{\cdot}$ denotes the floor function. This way, we write $X_i(t) = 1$ ($0$) if channel $i\in \cN$ is available (unavailable) for the SU with probability $p_i$ ($1-p_i$). 

In this setup, the rate at which the SU can transmit files through channel $i\in\cN$ at any time instant $t\geq0$, also termed as the \textit{instantaneous throughput} of the channel $i$, is given by $r_iX_i(t)$, with its \textit{throughput} \cite{Optimal_802.11,Srikant_2019} by $\E[r_i X_i(t)] = r_ip_i$. We denote by $i^* \triangleq \argmax_{k\in\cN} r_kp_k$ the channel with the \textit{maximum throughput}. For simplicity, we assume that this channel is unique, i.e., $r_{i^*}p_{i^*} \!>\! r_k p_k$ for all $k \in \cN \!\setminus\! \{ i^*\}$. In what follows, we introduce some basic notations and expressions regarding \textit{policies} under this OSA framework. 

\subsection{Policies for File Transfer}\label{subsection:policies for file transfer}

\begin{figure}[ht]
    \centering
    \includegraphics[width=0.85\columnwidth]{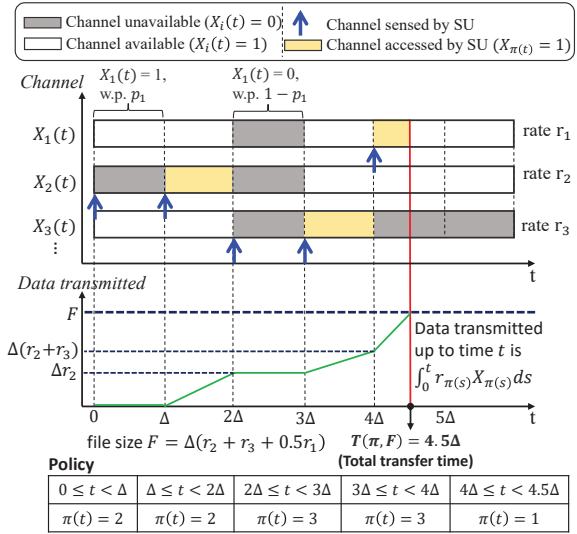} 
    \vspace{-2mm}
    \caption{File transfer via the OSA framework. The SU senses channels according to its policy, and accesses the channel if it is available. Upon gaining access, it begins transmitting data at the corresponding channel rate. Transmission ends (red line) as soon as the amount of data transmitted (green line) equals the file size $F$.}
    \label{fig:illustration}\vspace{-3mm}
\end{figure}

We define a \textit{policy} at time $t$ to be a mapping $\pi: \R_+ \mapsto \cN$ where $\pi(t)=i$ indicates that the SU has chosen channel $i$ to access during the time period $\left[\floor{\frac{t}{\Delta}}\!\Delta, \left(\floor{\frac{t}{\Delta}} \!\!+\!\!1 \right)\!\Delta\right)$. From our standing assumption, a policy therefore only changes at $t \!\in\! \{0,\Delta,2\Delta,\cdots\}$, and all policies ensure that the file transfer for any finite size $F$ will eventually be completed. This way, the policy $\pi(t)$ is a piecewise constant function (mapping), defined at all time $t \geq 0$. For a given policy $\pi$, let $T(\pi,F)$ denote the \emph{transfer time} of a file of size $F$ --- the entire duration of time to complete the file transfer, which is written as
\begin{equation}\label{stopping time}
    T(\pi,F) = \min_{T\geq 0}\left\{\int_0^{T} r_{\pi(t)} X_{\pi(t)}(t)dt ~~\geq F \right\}.
\end{equation}
Figure \ref{fig:illustration} explains the file transfer progress via OSA model.

The objective of our OSA framework is to minimize the expected transfer time $\E[T(\pi,F)]$ over the set of all policies $\pi$. Policies can be \textit{static}, where the SU only senses and transmits via one (pre-determined) channel $i$ throughout the file transfer, that is, $\pi(t) = i$ for all $t\geq0$. For such static policies, we denote by $T(i,F)$ their transfer time for file size $F$. The channel that provides the minimum expected transfer time is then called \textit{static optimal} given by
\begin{equation}\label{eqn:static_optimal_channel}
    i_{so}(F) = \text{arg min}_{k \in \cN} ~\E [T(k,F)], \quad \text{(static optimal)}
\end{equation}
\noindent and we denote $T(i_{so},F)$ the corresponding transfer time for this static optimal policy. Note that the static optimal channel $i_{so}(F)$ depends on the file size $F$ and can vary for different file sizes. Policies can also be \textit{dynamic}, in which an SU is allowed to change the channels it chooses to sense throughout the course of the file transfer. Given a file size $F$, the policy with the
minimum expected transfer time over the set of all policies $\Pi(F)$ is called the \textit{dynamic optimal} policy given by \begin{equation}\label{eqn:dynamic_optimal_channel}
    \pi^*(F) = \text{arg min}_{\pi \in \Pi(F)} ~\E [T(\pi,F)].~ \text{(dynamic optimal)}
\end{equation} 
Lastly, we define the \textit{max-throughput policy} as the static policy with the channel $i^*$, which maximizes the long-term throughput. In the next section, we take a closer look at the max-throughput policy and static policies in general.

\section{Static optimal policy}\label{static analysis}

Recent works in the OSA literature focus on estimating channel parameters $p_i$'s, with the goal of eventually converging to the policy $i^* \!=\!\! \argmax_{k\in\cN} \! r_kp_k$ which provides the maximum throughput \cite{Poor_2011,dai2012efficient,Optimal_802.11,Srikant_2018,Srikant_2019}.\footnote{While \cite{Optimal_802.11,Srikant_2018,Srikant_2019} deal with link rate selection problem to select best rate in one channel to maximize the expected throughput, the mathematical model of link rate selection problem is essentially the same as the standard OSA setting for choosing the max-throughput channel, as considered in our setting.} They focus on minimizing the `regret' in the MAB model, defined as the difference between the cumulative reward obtained by the online algorithm and the max-throughput policy (the optimal policy in hindsight).

The essential assumption behind all these approaches is that the SU always fully dedicates $\Delta$ seconds in each time interval for file transfer. Channel $i^*$ appears as a good candidate since it provides the largest expected data transfer $\Delta r_{i^*} p_{i^*}$ across every time interval. This is further supported by the well-known Wald's equation with the \textit{i.i.d} reward assumption at each time interval, suggesting that $\E[T(i,F)] \!=\! F/r_ip_i$ for each channel $i\!\in\!\cN$, which is then minimized by $i^*$. When policies are dynamic, however, the rewards are not identically distributed since the transfer rates of the dynamically accessed channels can be different, making Wald's equation inapplicable. Surprisingly, it is not applicable for static policies either. As typically is the case in delay-sensitive applications \cite{shiang2008queuing,wu2014learning,dimitriou2018stable, huang2019dynamic}, the file sizes are often not that large, rendering their transfer times small enough that an SU may not need to utilize the whole $\Delta$ seconds for data transfer in each time interval. The reward summands are still not identically distributed, causing Wald's equation to fail in general.

Our key observation in this paper is that choosing channel $i^*$ may not be the best option to minimize the expected transfer time. In this section, we limit ourselves to the set of static policies of the form shown in \eqref{eqn:static_optimal_channel} and analyze the resulting expected transfer time in the OSA network. We use this to compare the performance gap between the max-throughput policy and the static optimal policy, and show that for a reasonable choice of channel statistics and file sizes, the static optimal policy performs significantly better than the max-throughput policy. 
We derive a closed-form expression of the expected transfer time of a file of size $F$ in each fixed channel by the following proposition.

\begin{proposition}\label{prop:expected_time_static} 
Given a file of size $F$, the expected transfer time $\E[T(i,F)]$ of the static policy for channel $i \in \cN$ is
\begin{equation}\label{eqn:time_static_policy}
     \E[T(i,F)] = \Delta \left(k_i/p_i + \Char_{\{\alpha_i > 0\}}(1-p_i)/p_i + \alpha_i\right),
\end{equation}
where $ k_i \triangleq \floor{F/\Delta r_i} \in \Z_+$ and $\alpha_i \triangleq  F/ \Delta r_i - k_i \in [0,1)$.
\end{proposition}
\begin{IEEEproof}
See Appendix \ref{app:proof_prop_3.1}.
\end{IEEEproof}

We have from Proposition \ref{prop:expected_time_static} that $\E[T(i,F)]$, the expected transfer time for any file size $F$ under the static policy on channel $i \in \cN$, can be explicitly written in terms of file size $F$, time duration $\Delta$ and channel statistics $r_i$ and $p_i$ of the chosen channel $i$. Substituting $k_i = F/\Delta r_i - \alpha_i$ in \eqref{eqn:time_static_policy} gives 
\begin{equation}\label{eqn:static_lower_bound}
    \E[T(i,F)] \!=\! \frac{F}{r_ip_i} \!+\! \Delta \Char_{\{\alpha_i>0\}}  (1\!-\!\alpha_i) \frac{1\!-\!p_i}{p_i} \!\geq\! \frac{F}{r_ip_i}.
\end{equation}
The inequality in \eqref{eqn:static_lower_bound} shows the expected transfer time of any static policy is no smaller than that given by Wald's equation.

\begin{figure}[ht]
    \centering
    \vspace{-3mm}
    \includegraphics[width=0.85\columnwidth]{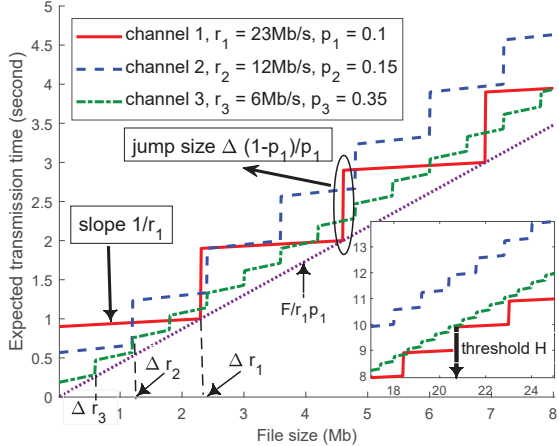}
    \vspace{-2mm}
    \caption{Expected transfer time from \eqref{eqn:time_static_policy} with duration $\Delta = 100$ ms (IEEE 802.22 standard). Channel $1$ has the maximum throughput. The purple dot-line $E[T] = F/r_1p_1$ corresponds to the Wald's equation for channel $1$. Downward arrow in the inset figure is the threshold $H$ in Proposition \ref{prop:threshold_optimal_policy}.}
    \label{fig:time_3channel_static}
    \vspace{-1mm}
\end{figure}

We use Figure \ref{fig:time_3channel_static} to illustrate the results in Proposition \ref{prop:expected_time_static}, where each line represents the expected transfer time via one channel over a range of file sizes from \eqref{eqn:time_static_policy}. We observe that the expected transfer time of channel $1$ (red line) is always above the Wald's equation of channel $1$ (purple dot-line). As shown in \eqref{eqn:time_static_policy}, the slope of each line (channel $i$) is $1/r_i$ and the `jump size' $\Delta(1-p_i)/p_i$ is equal to the expected waiting time till the channel is available. The jumps in the plot for each channel $i$, representing the waiting times, occur at exactly the instances where file size is an integer multiple of $\Delta r_i$, and come into play especially when there is still a small amount of remaining file to be transferred at the end of a $\Delta$ time interval.

By definition, the static optimal policy provides the minimum expected transfer time over all static policies including the max-throughput policy itself. While it is true for all file sizes, in some cases with certain file sizes, these two policies may coincide.

\begin{proposition}
\label{prop:threshold_optimal_policy}
The max-throughput policy coincides with the static optimal policy, that is, $i_{so}(F) = i^*$, for any file size $F$ satisfying at least one of the two conditions below:
\begin{enumerate}
    \item $F$ exceeds a threshold $H$, where
\begin{equation} \label{eqn:static_policy_threshold}
    H = \frac{\Delta (1-p_{i^*})/p_{i^*}}{1/r_hp_h - 1/r_{i^*}p_{i^*}},
\end{equation}
and $h = \argmax_{j \in \cN \setminus \{i^*\}} r_jp_j$ is the channel with the second largest throughput.
    \item $F$ is an integer multiple of $\Delta r_i^*$, i.e., $F \!=\! k \Delta r_{i^*}$ for some $k \!\in\! \Z_+$. 
\end{enumerate}
\end{proposition}
\begin{IEEEproof}
See Appendix \ref{app:proof_prop_3.2}.
\end{IEEEproof} 

Outside of Proposition \ref{prop:threshold_optimal_policy}, however, there are many instances where the max-throughput channel is not static optimal and other channels can perform better for smaller file sizes. In such cases, we would like to discuss how much time the static optimal policy can save against the max-throughput policy.

\begin{corollary}\label{corollary:3.3}
Let $m_i \triangleq p_i/p_{i^*}$ for $i\in \cN \setminus \{i^*\}$. Consider a file of size $F \in (k\Delta r_{i^*}, (k+1)\Delta r_{i^*})$ for some $k \in \N$. Then, we have
\begin{equation}\label{eqn:lower_bound_static_general}
\begin{split}
   \frac{\E[T(i_{so},F)]}{\E[T(i^*,F)]} 
   \leq \!\! \min_{\substack{i\in \cN \setminus \{i^*\}}}\! \left\{1,\frac{F/\Delta r_ip_i + (1-p_i)/p_i}{(k+1)(m_i/p_i-1)}\right\}.
\end{split}
\end{equation}
\end{corollary}
\begin{IEEEproof}
See Appendix \ref{app:proof_coro_3.3}.
\end{IEEEproof}

Note that by definition $m_i/p_i = 1/p_{i^*} > 1$ so that the upper bound on the ratio $\E[T(i_{so},F)]/\E[T(i^*,F)]$ in \eqref{eqn:lower_bound_static_general} is always in the interval $(0,1]$. Moreover, smaller ratio means better performance of the static optimal policy against the max-throughput policy. To gauge how the parameters of the max-throughput channel could affect the performance of the static optimal policy, suppose we fix $F,\Delta,r_i,p_i$ for all $i\!\in\! \cN\! \setminus\! \{i^*\}$ and the maximum throughput $r_{i^*} p_{i^*}$, while treating $p_{i^*}$ as a variable. The upper bound in \eqref{eqn:lower_bound_static_general} is then monotonically decreasing in $m_i = p_i/p_{i^*}$, and can even approach to $0$ if at least one of $m_i$ is really large, resulting in the huge performance gain of the static optimal channel compared to that of the max-throughput channel. This implies that accessing channel $i^*$ can take much longer time to transmit a file than other channels if its available probability $p_{i^*}$ is very small, which is common in outdoor networks where the max-throughput channel $i^*$ has very high rate but with low available probability  \cite{bicket2005bit}.

Our static optimal policy shows better performance against the max-throughput policy for small files and small $p_{i^*}$. Since the static optimal channel depends on the file size, choosing channels \emph{dynamically} according to its remaining file size can further reduce the expected transfer time. We next formulate the file transfer problem as an instance of the stochastic shortest path (SSP) problem and analyze the performance of the dynamic optimal policy. 

\section{Dynamic Optimal Policy}\label{mdp formulation and dynamic policy}

Now that we have analyzed the static policies, we turn our attention to feasible dynamic policies for our file transfer problem. We start by first formulating the file transfer problem as a stochastic shortest path (SSP) problem, in which the agent acts dynamically according to the stochastic environment to reach the predefined destination as soon as possible.
Then, we translate this SSP problem into an equivalent shortest path problem, which helps us derive the closed-form expression of the expected transfer time for any given dynamic policy, and we utilize this to obtain the performance analysis of the dynamic optimal policy against the max-throughput policy.

\subsection{Stochastic Shortest Path Formulation}\label{MDP setting}
The SSP problem is a special case of the infinite horizon Markov decision process \cite{bertsekas2019reinforcement}. To make this section self-contained, we explain our problem as a SSP problem.

\textit{State Space} and \textit{Action Space:} We define the state $s \in \cS \triangleq \R_+$ of our SSP as the remaining file size yet to be transmitted. The action $i \in \cN$ is the channel chosen to be sensed at the beginning of each time interval. The objective of our problem is to take the optimal action at each state $s$ which minimizes the expected time to transmit the file of size $F$.

\textit{State Transition:} Denote by $P_{s,s'}(i)$ the transition probability that the SU moves to state $s'$ after taking action $i$ at state $s$. From any given state $s \in \cS \!\setminus\! \{0\}$, the next state under any action $i\in \cN$ depends on the availability of channel $i$. Since the channel is available or unavailable according to an \textit{i.i.d} (over time) Bernoulli distribution, the next state is either the same as the current one if channel $i$ is unavailable, i.e. $P_{s,s}(i) = 1-p_i$; or the next state is $(s\!-\!\Delta r_i)^+ \triangleq \max\{0,s\!-\!\Delta r_i\}$ if channel $i$ is available, i.e. $P_{s,(s\!-\!\Delta r_i)^+}(i) = p_i$. State $0$ is a termination state since there is no file transmission remaining.

\textit{Cost Function:} The cost $c(s,i,s')$ is the amount of time spent in transition from state $s$ to $s'$ after sensing channel $i$. Since the SU can only sense channels at intervals of size $\Delta$, sensing an unavailable channel costs a $\Delta$ second waiting period until the SU can sense next, that is, $c(s,i,s) = \Delta$ for all $s\in \cS \!\setminus\! \{0\},i \in \cN$. Similarly, if the sensed channel is available, the time spent in transmitting is also $\Delta$ seconds, unless the SU finishes transmitting the file early. In the latter case the cost of transmission is $c(s,i,0) = s/r_i$. Overall, the cost of a successful transmission can be written as $c(s,i,(s\!-\!\Delta r_i)^+) = \min\{\Delta,s/r_i\}$ for all $s\in \cS \!\setminus\! \{0\},i \in \cN$. Once the remaining file size reduces to $0$, the SU will end this file transmission immediately with no additional cost incurred, so that $c(0,i,0) \!=\! 0$ for any $i \in \cN$.

Table \ref{tab:ssp_setting} summarizes the state transition and cost function for our file transfer problem. All other cases except the two cases in Table \ref{tab:ssp_setting} have zero transition probability and zero cost.
\begin{table}[!ht]\vspace{-4mm}  
    \centering
    \caption{SSP setting of our file transfer problem}
    \label{tab:ssp_setting}
    \begin{tabular}{|c|c|c|c|c|c|}
         \hline
         current state & action& next state & transition & cost \\
         \hline
         $s>0$ & $i$ & $s$ & $1-p_i$ & $\Delta$ \\
         \hline
         $s>0$ & $i$ & $(s\!-\!\Delta r_i)^+$ & $p_i$ & $\min\{\Delta,s/r_i\}$ \\
         \hline
    \end{tabular}
    \vspace{-2mm}
\end{table}

Our dynamic policy\footnote{There always exists an optimal policy $\pi^*$ to be deterministic in the SSP problem, as proved in Proposition 4.2.4 \cite{bertsekas2019reinforcement}. Thus, we restrict ourselves to the class of deterministic policies in this paper.} is written as a mapping $\pi : S \to \cN$, where $\pi(s)\in\cN$ denotes the channel chosen for sensing when the current state (remaining file size) is $s$. For any policy $\pi$, we have $T(\pi,0) = 0$ at the termination state. Our goal in this SSP problem is to find the dynamic optimal policy $\pi^*(F)$ that minimizes the expected transfer time for the file size $F$, which can be derived from a variety of methods such as value iteration, policy iteration and dynamic programming \cite{bertsekas2019reinforcement}.

\subsection{Performance Analysis}\label{performance analysis}
\begin{figure}[H]
    \vspace{-0mm}
    \subfloat{\includegraphics[width=\columnwidth]{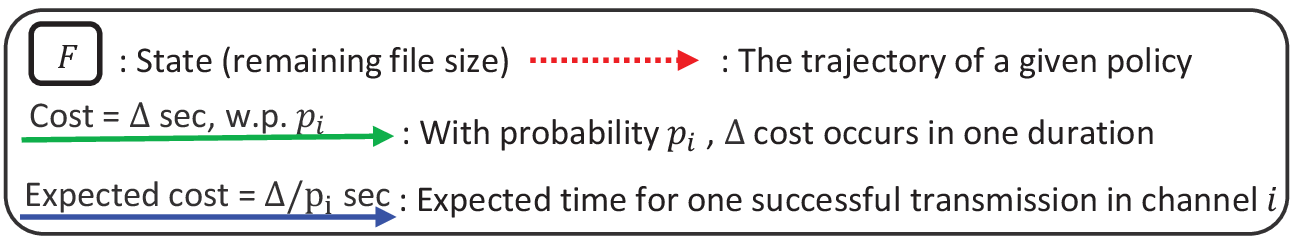}}\qquad
    \subfloat[(a) Transformation of the SSP problem.\label{fig:deterministic_figure}]{\includegraphics[,width=0.3\columnwidth]{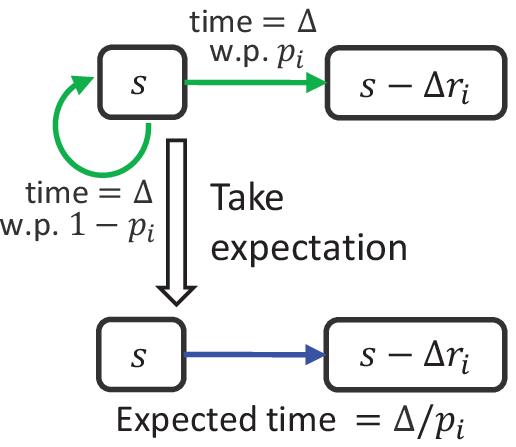}}\qquad
    \subfloat[(b) Equivalent shortest path diagram.\label{fig:mdp_instance}]{\includegraphics[width=0.63\columnwidth]{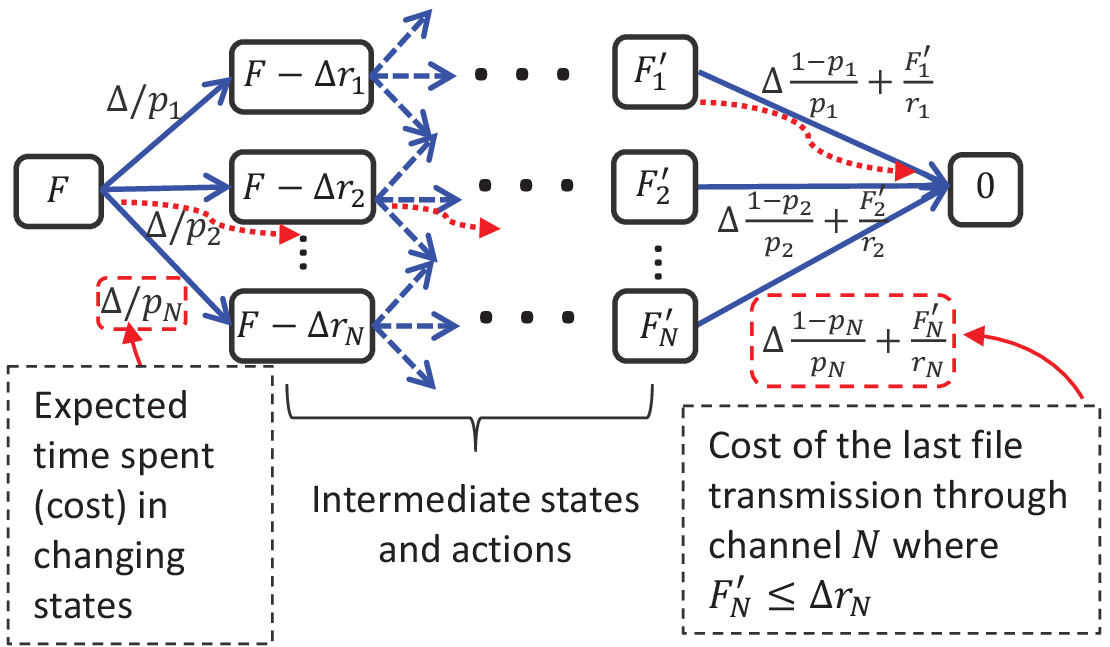}}
     \caption{Illustration to translate the file transfer problem into an equivalent shortest path problem.}
     \label{fig:equivalent_shortest_path}\vspace{-2mm}
\end{figure}

For ease of exposition, we introduce additional notations here. By a \emph{successful transmission}, we refer to state transitions of the form $s \to s'$. This is denoted by the horizontal green line in Figure \ref{fig:equivalent_shortest_path}(a) connecting states $s$ and $s' \triangleq s - \Delta r_i > 0$, and should be distinguished from the self-loop $s \to s$, which implies the sensed channel was unavailable. As shown in Figure \ref{fig:equivalent_shortest_path}(a), taking expectation helps get rid of these self-loops by casting the original SSP to a deterministic shortest path problem in expectation. The cost associated with each link is then the expected time it takes to transit between the states. Figure \ref{fig:equivalent_shortest_path}(b) shows the underlying network for the shortest path problem, where each link is a channel chosen to be sensed and each path from \emph{source} $F$ to \emph{destination} $0$ corresponds to a policy $\pi \in \Pi(F)$. The path-length or the number of links traversed from $F$ to $0$ under any given policy $\pi$ then becomes the total number of successful transmissions needed by that policy to complete the file transfer, which we denote by $|\pi|$.

For any policy $\pi \in \Pi(F)$ and $n\in\{1,\cdots,|\pi|\}$, let $F_n$ denote the remaining file size right before the $n$-th successful transmission. Then for all $n\in \{2, \cdots, |\pi|\}$, we have the recursive relationship: $F_n = F_{n-1} - \Delta r_{\pi(F_{n-1})}$, starting with $F_1 \!=\! F$ and ending with $F_{|\pi|+1}=0$. Given a file size $F$, each policy $\pi \in \Pi(F)$ can then be written in a vector form as $\pibf = [\pi(F_1), \pi(F_2), \dots \pi(F_{|\pi|})]^T$.
With this in mind, we can derive a closed-form expression of the expected transfer time for any dynamic policy in the following proposition.

\begin{proposition}\label{prop:expected_time_dynamic_policy}
    Given a file of size $F$, the expected transfer time of a dynamic policy $\pi$ is written as
\begin{equation}\label{eqn:recursion_closed_form}
\begin{split}
       \E[T(\pi,F)]  \!=\! \Delta \!\!\left(\!\sum_{n=1}^{|\pi|\!-\!1}\! \frac{1}{p_{\pi(F_{n})}}\!\!\right)\! \!+\!\Delta\frac{1\!-\!p_{\pi(F_{|\pi|})}}{p_{\pi(F_{|\pi|})}}\!+\! \frac{F_{|\pi|}}{r_{\pi(F_{|\pi|})}}. 
\end{split}
\end{equation}
\end{proposition}
\begin{IEEEproof}
See Appendix \ref{app:proof_prop_4.1}.
\end{IEEEproof}

In \eqref{eqn:recursion_closed_form}, the first summation is the cumulative expected transmission time, or the cost, to the $(|\pi|-1)$-th successful transmission, with the last two terms being the expected transmission time of the last successful transmission. Proposition \ref{prop:expected_time_dynamic_policy} also includes the expected transfer time of the static policy as a special case. Recall that $k_{i} = \floor{F/\Delta r_i}$ and $\alpha_{i} = F/\Delta r_i - k_i$ in Proposition \ref{prop:expected_time_static}. When applied to a static policy for any channel $i \in \cN$, we have $|\pi| = k_i + \Char_{\{\alpha_i > 0\}}$ and $p_{\pi(F_n)} = p_i$ for all $n = 1,2,\cdots,|\pi|$. The recursive relationship becomes: $F_n = F_{n-1} - \Delta r_i$, implying that $F_n = F - (n-1)\Delta r_i$. Then, we have $F_{|\pi|} = F - (|\pi|-1)\Delta r_i = \alpha_i\Delta r_i$ if $\alpha_i > 0$. Otherwise, $F_{|\pi|} = \Delta r_i$. Substituting these into \eqref{eqn:recursion_closed_form} gets us \eqref{eqn:time_static_policy}.

The common folklore around the max-throughput policy is that it would lead to the minimal file transfer time of $F/r_{i^*}p_{i^*}$. Our next result shows this is too optimistic and not achieved in general even under the dynamic optimal policy. 

\begin{proposition}\label{prop:lower_bound_coincide}
For any file size $F$ and any dynamic policy $\pi \in \Pi(F)$, we have $\E[T(\pi,F)] \geq \E[T(\pi^*,F)] \geq F/r_{i^*}p_{i^*}$.
Moreover, $i^* = \pi^*(F)$ for $F = k\Delta r_{i^*}, k \in \Z_+$.
\end{proposition}
\begin{IEEEproof}
See Appendix \ref{app:proof_lemma_4.2}.
\end{IEEEproof}

As shown in \eqref{eqn:static_lower_bound}, $F/r_{i^*}p_{i^*}$ is always the lower bound on the transfer time for any static policy. Proposition \ref{prop:lower_bound_coincide} strengthens this by showing that the same is true even for the dynamic optimal policy. Similar to condition (b) in Proposition \ref{prop:threshold_optimal_policy} for the static optimal policy, the dynamic optimal policy $\pi^*(F)$ also coincides with the max-throughput policy $i^*$ when the file size is an integer multiple of $\Delta r_{i^*}$, while we no longer have the finite threshold $H$ as in Proposition \ref{prop:threshold_optimal_policy}(a). We next give bounds to quantify the performance of the dynamic optimal policy with respect to the max-throughput policy. 
\begin{corollary}\label{corollary:dynamic_lower_upper_bound}
Let $m_i \triangleq p_i/p_{i^*}$ for $i\in \cN \setminus \{i^*\}$. Consider a file of size $F \in (k\Delta r_{i^*}, (k+1)\Delta r_{i^*})$ for some $k \in \N$. Then, we have
\begin{equation*}
\begin{split}
    &1/(1+\Delta \Char_{\{\alpha_i > 0\}}(1-p_{i^*})r_{i^*}/F)\leq \frac{\E[T(\pi^*,F)]}{\E[T(i^*,F)]}  \\
    &\leq \!\!\!\!\min_{\substack{i \in \cN \!\setminus\! \{i^*\}}}\!\!\!\left\{\!\!1,\!\frac{F/\Delta r_ip_i \!+\! (1\!-\!p_i)/p_i \!-\! k m_i(r_{i^*}p_{i^*}\!\!-\!r_ip_i)/r_ip_i^2}{(k+1)(m_i/p_i\!-\!1)}\!\right\}\!\!.
\end{split}
\end{equation*}
\end{corollary}
\begin{IEEEproof}
See Appendix \ref{app:proof_coro_4.4}.
\end{IEEEproof}

To better understand Corollary \ref{corollary:dynamic_lower_upper_bound} we analyze how the parameters of the max-throughput channel could impact the performance of the dynamic optimal policy. Similar to Corollary \ref{corollary:3.3}, small value of ${\E[T(\pi^*,F)]}/{\E[T(i^*,F)]}$ implies that the dynamic optimal policy offers significant saving in time over the max-throughput policy. We note that Corollary \ref{corollary:dynamic_lower_upper_bound} tightens the upper bound with an extra negative term in the numerator, compared to that in Corollary \ref{corollary:3.3}, potentially providing greater savings in time as we extend the policy from static optimal to the dynamic optimal.

In contrast to Proposition \ref{prop:threshold_optimal_policy} that max-throughput policy is good enough for $F \geq H$, Corollary \ref{corollary:dynamic_lower_upper_bound} tells us that there is always some reduction in file transfer time even for large file size $F$ under the dynamic optimal policy. This is because the extra negative term in the numerator can be large, since $k=\floor{F/\Delta r_{i^*}}$ could be big for large $F$, implying that the second argument in the $\min\{\cdot,\cdot\}$ function may no longer be increasing in $F$. Note however that the reduction in transfer time would be minimal for large file sizes since the lower bound in Corollary \ref{corollary:dynamic_lower_upper_bound} will rise to $1$ as $F$ goes to infinity.

\section{Practical considerations}\label{online learning}
While the dynamic optimal policy gives a smaller expected transfer time, we face a scaling problem when the file size can differ from each, effectively changing the underlying `graph' in the corresponding shortest path problem. This warrants re-computation of the dynamic optimal policy for each file size, which would be unacceptable in reality. In this section, we discuss the policy re-usability issue and propose a heuristic policy balancing the performance and the computational cost. Then, we consider the case where the SU has no information about the channel parameters beforehand and it must sense and access channels \textit{on the fly} in order to find the optimal policy, thereby extending the problem to an \textit{online} setting. Finally, we will present an mixed-integer programming formulation for the dynamic optimal policy in our file transfer problem that could reduce the computational cost.

\subsection{Performance-Complexity Trade Off}\label{subsection: reusable}
Transmitting different-sized files is very common in the real world. For example, a short text-only email only takes up $5$KB, one five-page paper is around $100$KB and the average size of webpage is $2$MB, all implying that the file size may vary greatly \cite{mantuano_2016}. However, due to the nature of the shortest path problem, a change in file size induces a change in the underlying graph. If the goal is to always determine the best solution, the only option is to recompute the dynamic optimal policy for every different file size. This would not be scalable in applications where minimal computation is required, and policies that can be promptly modified and reused across different file sizes with performance guarantees would be highly desirable.

To avoid heavy computation for each file size to obtain the dynamic optimal policy, we propose a heuristic policy that utilizes the max-throughput policy and the static optimal policy in order to reduce the computational cost, while still maintaining considerable performance gain. Note that the max-throughput policy coincides with the dynamic optimal policy when the file size is an integer multiple of $\Delta r_{i^*}$ according to Proposition \ref{prop:lower_bound_coincide}. We also know that the static optimal policy significantly outperforms the max-throughput policy especially for smaller file sizes. Combining these two policies, by transmitting file through max-throughput channel until the remaining file size becomes `small'\footnote{We can choose remaining file size to be smaller than $\Delta r_{i^*}$ to apply the static optimal policy.} so as to apply the static optimal policy for the rest, will strike the right balance between computational complexity and achievable performance gain. From the computational-cost-saving perspective, max-throughput policy is fixed and known to the SU. The closed-form expression $E[T(i,F)]$ for static policy (channel $i$) is also known to the SU, which means the computational cost of the static optimal policy by taking the minimum delay over all $N$ channels is much smaller than that of the dynamic optimal policy. 

With this motivation in mind, we divide file size $F := F_1 + F_2$ into two parts: $F_1 = n \Delta r_{i^*}$ ($n = 0,1,\cdots,\floor{F/\Delta r_{i^*}}$) and $F_2 = F- F_1$. The heuristic policy $\pi_{heu}$ is defined as follows: The SU first transmits the file of size $F_1$ through the max-throughput channel $i^*$ and then sticks to the static optimal policy $i_{so}(F_2)$ for the remaining file of size $F_2$.\footnote{For $F$ being integer multiple of $\Delta r_{i^*}$, we have $F_1 = F$ and $F_2 = 0$, then the heuristic policy coincides with the max-throughput policy, which is also the dynamic optimal policy in view of Proposition \ref{prop:lower_bound_coincide}.} We show in the Appendix \ref{app:proof_coro_4.4} that the upper bound of ratio $E[T(\pi_{heu},F)]/E[T(i^*,F)]$ is monotonically decreasing in $n$. Therefore, $n = k$ gives us the smallest upper bound of the ratio (the same upper bound in Corollary \ref{corollary:dynamic_lower_upper_bound}). Moreover, we have explained after Corollary \ref{corollary:dynamic_lower_upper_bound} that the upper bound is smaller than that of the static optimal policy in Corollary \ref{corollary:3.3}. These arguments suggest that the heuristic policy $\pi_{heu}$ with $n=k$ can potentially offer smaller delay than other candidates with different values of $n$. Besides, our heuristic policy can further reduce the computational cost for a set of files sharing the same remaining file size $F_2$, for which the static optimal policy $i_{so}(F_2)$ has already been found and no further recomputation is needed for this set of files.

\subsection{Unknown Channel Environment}\label{subsection:online_implementation}
We now consider the setting where the SU does not know the available probability $p_i$ for any channel $i\in\cN$ and only knows the rate $r_i$ --- a commonly analysed setting in the OSA literature \cite{Poor_2011,dai2012efficient}. The SU has no alternative but to observe the states of these channels when it tries to access them, and build its own estimations of channel probabilities. In this extended setting, we study our problem as an online shortest path problem, which has been widely studied in \cite{gai2012combinatorial,chen2013combinatorial,talebi2017stochastic} for different kinds of cost functions. \cite{talebi2017stochastic} proposed a Kullback-Leibler source routing (\textit{KL-SR}) algorithm to an online routing problem with geometrically distributed delay in each link, which coincides with our link cost in the underlying graph of the shortest path problem shown in Figure \ref{fig:equivalent_shortest_path}(b).

For our purpose, we modify \textit{KL-SR} algorithm; key differences being that we let the file size vary across the episodes, allowing a different underlying graph of the shortest path problem for each episode instead of the fixed underlying graph of the shortest path problem in \cite{talebi2017stochastic}. Algorithm \ref{alg:online_algorithm} describes our online implementation,
where $F^k$ is the file size to be transferred in the $k$-th episode. $n_i(k)$ denotes the number of times channel $i$ has been sensed before the $k$-th episode and $\bar{p}_i(k)$ is the empirical average of channel $i$'s available probability throughout the $k\!-\!1$ episodes so far. With $n_i(k)$ and $\bar{p}_i(k)$, the estimated available probability $\hat{p}_i(k)$ of channel $i$ is then derived from the KL-based index in \cite{talebi2017stochastic}. As mentioned in line $1$ in Algorithm \ref{alg:online_algorithm}, the SU can choose one of the various policies according to which it wishes to perform the file transfer, i.e., dynamic optimal policy $\pi^*$, static optimal policy $i_{so}$, max throughput policy $i^*$ or the heuristic policy $\pi_{heu}$, and then stick to that policy.
Let $\tilde{\E}[T(\pi,F^k)]$ be the estimated expected transfer time of policy $\pi$ at the $k$-th file by using the estimated parameter $\hat{p}_i(k)$ instead of $p_i$ for all $i \in \cN$ in \eqref{eqn:recursion_closed_form}. In line $7$ in Algorithm \ref{alg:online_algorithm}, $\pi^k$ will be computed as $\pi^k \!=\! \argmin_{\pi \in \Pi(F^k)}\tilde{\E}[T(\pi,F^k)]$ for dynamic optimal policy; $\pi^k \!=\! \argmin_{i \in \cN}\tilde{\E}[T(i,F^k)]$ for static optimal policy and $\pi^k \!=\! \argmax r_i \hat{p}_i(k)$ for max-throughput policy. For heuristic policy $\pi_{heu}$,\! $\pi^k\!$ will be computed in the same way as described in Section \ref{subsection: reusable} with estimated parameters $\{\hat{p}_i(k)\}_{i \in \cN}$.
\vspace{-1mm}
\begin{algorithm}[H]
	\caption{Online file transfer algorithm}
	\label{alg:online_algorithm}
	\begin{algorithmic}[1]
	    \State Choose the type of policy to use: $\pi^*$, $i_{so}$, $i^*$ or $\pi_{heu}$.
	    \State Apply static policy $i$ in the $i$-th episode to transmit the file of size $F^i$ for $i=1,2,\cdots,N$, and update $n_i(N+1)$, $\bar{p}_i(N+1)$ for $i \in \cN$ at the end of the $N$-th episode 
	    \For {$k \geq N+1$}
	        \State Compute the estimated channel statistics $\{\hat{p}_i(k)\}_{i\in \cN}$ according to the KL-based index in \cite{talebi2017stochastic}.
	        \State Given a file of size $F^k$, compute the policy $\pi^k$ with $\{\hat{p}_i(k)\}_{i \in \cN}$ and observe channel status for the whole file transfer process in this episode. 
	        \State Update $n_i(k+1)$, $\bar{p}_i(k+1)$ for $i \in \cN$.
	    \EndFor
	\end{algorithmic}
\end{algorithm}\vspace{-2mm}

The performance of Algorithm \ref{alg:online_algorithm} with varying file sizes (assuming bounded file size) is measured by its \textit{regret} $\E[R(K)]$, which is defined as the cumulative difference of expected transfer time between policy $\pi^k$ at $k$-th file and the targeted optimal policy $\pi_{tar} \in \{\pi^*, i_{so}, i^*, \pi_{heu}\}$ up to the $K$-th file. The regret analysis is nearly the same as Theorem $5.4$ in \cite{talebi2017stochastic} and the regret bound of Algorithm \ref{alg:online_algorithm} is given in the following theorem.
\begin{theorem}\label{thm:regret_bound}
The gap-dependent regret bound under Algorithm \ref{alg:online_algorithm} is
\begin{equation}\label{eqn:regret_bound}
    \E[R(K)] \leq \frac{360NH f(K)}{\Delta_{\min}p_{\min}^2} + 2D\left(4H + \sum_{i=1}^n \frac{1}{\epsilon^2 p^2_i}\right),
\end{equation}
where $f(K) = \log(K)+4\log(\log(K))$, $H = F_{\max}/r_{\min}$, $F_{\max}$ is the largest file size, $r_{\min} \!=\! \min_{i \in [N]} r_i$, $p_{\min} \!=\! \min_{i \in [N]} p_i$, $D = \max_{\pi} \E[T(\pi,F_{\max})]$ is the longest expected transfer time, $\epsilon = (1-2^{-\frac{1}{4}})\Delta_{\min}/D$ and $$\Delta_{\min} \!=\! \min_{\{F \in (0,F_{\max}],\pi\in \Pi(F)/\pi_{tar}\}} \E[T(\pi,F)] \!-\! \E[T(\pi_{tar}(F),F)]$$ is the smallest non-zero difference of expected transfer time between any sub-optimal policy $\pi$ and the targeted optimal policy $\pi_{tar}$.
\end{theorem}
\begin{IEEEproof}
See Appendix \ref{app:proof_thm_5.1}.
\end{IEEEproof}
The regret \eqref{eqn:regret_bound} scales linearly with the number of channels $N$, instead of the number of edges in the online shortest path problem \cite{talebi2017stochastic}, because each edge in our setting (see Figure \ref{fig:equivalent_shortest_path}) is chosen from one of $N$ channels while each edge in \cite{talebi2017stochastic} is treated as a different `arm'.

\section{Numerical Results}\label{simulation}
In this section, we present numerical results for file transfer time under four different policies in online settings (unknown $p_i$'s), using three different channel scenarios as in \cite{bicket2005bit,Optimal_802.11,Srikant_2019}. Through these results, we show the significant time reduction achieved by the dynamic optimal, static optimal and heuristic polices over the max-throughput channel, in line with theoretical analysis.

We consider the experimental setup as an IEEE 802.22 system with $8$ different channels. The time duration $\Delta$ is set to $100$ ms, per IEEE 802.22 standard \cite{802.22_standard}. We use three different channel scenarios: \textit{gradual, steep} and \textit{lossy} \cite{bicket2005bit,Optimal_802.11,Srikant_2019}. \textit{Gradual} refers to a case where the available probability of the max-throughput channel is larger than $0.5$. \emph{Steep} is characterized by the available probability of each channel being either very high or very low. \textit{Lossy} means that the available probability of the max-throughput channel is smaller than $0.5$. The channel parameters in the above three channel scenarios are given in Table \ref{tab:parameter}. 
\begin{table}[!ht]\vspace{-4mm}
    \centering
    \caption{Channel parameters in three channel scenarios}
     \label{tab:parameter}
    \begin{adjustbox}{width=\columnwidth}
    \begin{tabular}{|c|c|c|c|c|c|c|c|c|}
         \hline
         channel $i$ & $1$ & $2$ & $3$ & $4$ & $5$ & $6$ & $7$ & $8$\\
         \hline
         $r_i$ (Mbps) & $1.5$ & $4.5$ & $6$ & $9$ & $12$ & $18$ & $20$ & $23$ \\
         \hline
         $p_i$ (\emph{gradual}) & $0.95$ & $0.85$ & $0.75$ & $0.65$ & $0.4$ & $0.3$ & $0.2$& $0.1$ \\
         \hline
         $p_i$ (\emph{steep}) & $0.9$ & $0.25$ & $0.2$ & $0.18$ & $0.17$ & $0.16$ & $0.15$ & $0.14$ \\
         \hline
         $p_i$ (\emph{lossy}) & $0.9$& $0.8$ & $0.7$ & $0.4$ & $0.3$ & $0.25$ & $0.2$ & $0.1$ \\
         \hline             
    \end{tabular}
    \end{adjustbox}
    \vspace{-1mm}
\end{table}


In our setting, the file size need not be fixed. Since larger file sizes naturally take more time to transmit, it makes sense to normalize our performance metric across the range of file sizes. We define our metrics as \emph{average time ratio} and \emph{average throughput}. For an arbitrary sequence of files $\{F^k\}_{k \in \Z_+}$, the average time ratio at the $K$-th episode is defined as $\frac{1}{K}\sum_{k=1}^K{T(\pi^k,F^k)}/{\E[T(i^*,F^k)]}$ and the average throughput is represented as $\frac{1}{K}\sum_{k=1}^K F^k/T(\pi^k,F^k)$. Here, $T(\pi^k,F^k)$ is the measured transfer time of a file of size $F^k$ applying the policy $\pi^k$ at the $k$-th episode. Policy $\pi^k$ is based on the estimated parameter, which is updated by the SU on the fly, as described in Section \ref{subsection:online_implementation}.

\begin{figure}[!ht]\vspace{-2mm}
    \centering
    \includegraphics[width=0.88\columnwidth]{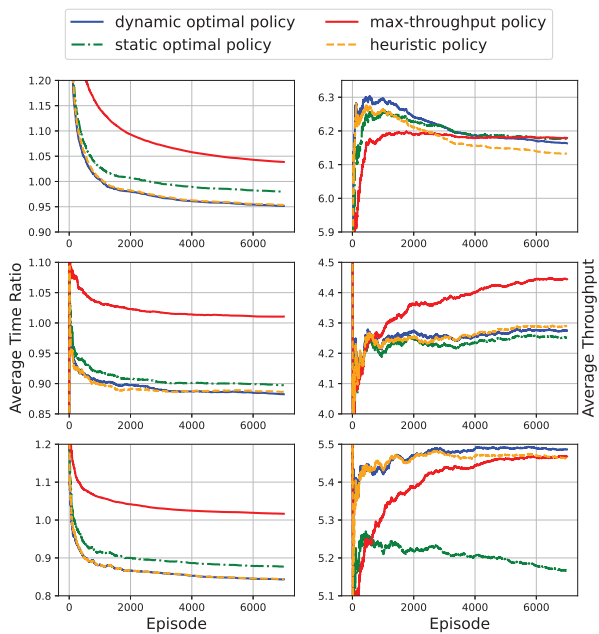}
    \vspace{-2mm}
    \caption{Average time ratio (left column) and average throughput (right column). Channel scenarios from top to the bottom: Gradual; Steep; Lossy.}
    \label{fig:online_simulation_gain}\vspace{-0mm}
\end{figure}

In our simulation, we generate $7000$ files from $(0,7]$ (Mb) uniformly at random to be used in Algorithm \ref{alg:online_algorithm}. The simulation is repeated $200$ times. 
We observe in Figure \ref{fig:online_simulation_gain} that the max-throughput policy achieves the largest average throughput while, counter-intuitively, has the longest transfer time in all channel cases except the lossy case. The reason is that the max-throughput policy computed by the SU is affected by the estimated parameters and can be the inferior policy, resulting in lower average throughput initially. For the lossy case on the right column in Figure \ref{fig:online_simulation_gain}, even though the red curve (max-throughput policy) is below the blue one for now, which is an effect of imperfect knowledge of channel parameters, we infer that the red curve will eventually exceeds all other curves as the SU will perfectly learn all the channel parameters.

Next we focus on the average time ratio in the left column of Figure \ref{fig:online_simulation_gain}. We first observe that all curves eventually flatten out, signifying the convergence of Algorithm \ref{alg:online_algorithm}. In the gradual case, the average time ratio is above $95\%$ for all three policies, implying that they don't obtain much reduction in time and the max-throughput channel is good to access when it is available for most of the time.
However, as shown in the \textit{steep} and \textit{lossy} cases respectively, the dynamic optimal policy and heuristic policy, as well as the static optimal policy, can save over $10\%$ time on average over the baseline. This observation is in line with Corollary \ref{corollary:3.3} and Corollary \ref{corollary:dynamic_lower_upper_bound} since the available probabilities of the max-throughput channel are very small in \textit{steep} and \textit{lossy} cases. Furthermore, the heuristic policy, in addition to keeping the complexity low, achieves similar transfer time to that of the dynamic optimal policy; at the same time performing better than the static optimal policy, as expected from Section \ref{subsection: reusable}.

\section{Conclusion and Future Work}
In this paper, we have developed a theoretical framework for the file transfer problem, where channels are modeled as independent Bernoulli process, to provide the accurate file transfer time for both static and dynamic policies. We pointed out that the max-throughput channel does not always minimize the file transfer time and provided static optimal and dynamic optimal policies to reduce the file transfer time. Throughout our analysis, we demonstrated that our approaches can obtain significant reduction in file transfer time over the max-throughput policy for small file sizes or when the max-throughput channel has very high rate but with low available probability, as typically the case in reality. Our future works include the extension to heterogeneous and \emph{Markovian} channels for minimal file transfer time, for which our SSP formulation for online learning scenario becomes no longer applicable.

\bibliographystyle{IEEEtran}
\bibliography{ref}

\appendices
\section{Proof of Proposition \ref{prop:expected_time_static}}\label{app:proof_prop_3.1}
Observe that any file size $F$ transmitted in channel $i$ can be written as
\begin{equation}\label{eqn:file_size_static}
    F = k_i \Delta r_i + \alpha_i \Delta r_i
\end{equation}
with $k_i$ being the number of intervals fully utilized for successful transmission, and $\alpha_i$ being the fraction of the $\Delta$ second interval utilized for file transfer toward the end. After choosing channel $i$, the SU first spends a random amount of time, denoted by $T_{wait}^n$ ($n = 1,2,\cdots,k_i$), waiting for channel $i$ to become available and starts the $n$-th transmission in that channel for $\Delta$ seconds. If the remaining portion $\Delta \alpha_i$ is not zero, the SU needs additional random waiting time $T_{wait}^{k_i+1}$ to complete the transfer. These random variables $\{T_{wait}^n\}$ are geometrically distributed and i.i.d over $n=1,2,\cdots,k_i+ \Char_{\{\alpha_i>0\}}$ with mean $E[T_{wait}^n] = \Delta (1-p_i)/p_i$. Let the constant $T_{tran} \triangleq F/r_i = \Delta(k_i+\alpha_i)$ be the total successful transmission time. Then the transfer time can be written as 
\begin{equation*}
    \begin{split}
        T &= \sum_{n=1}^{k_i} \left(T_{wait}^n + \Delta\right) + \Char_{\{\alpha_i > 0\}} T_{wait}^{k_i+1} + \Delta \alpha_i \\
        &= T_{tran} + \!\!\sum_{n=1}^{k_i+ \Char_{\{\alpha_i > 0\}}}\!\!  T_{wait}^n.
    \end{split}
\end{equation*}
Taking the expectation of the equation above yields \eqref{eqn:time_static_policy}.

\section{Proof of Proposition \ref{prop:threshold_optimal_policy}}\label{app:proof_prop_3.2}
From $\Char_{\{\alpha_i>0\}}(1-\alpha_i) \in [0,1]$ and \eqref{eqn:static_lower_bound}, we have the upper bound and the lower bound of $E[T(i,F)]$ as follows:
\begin{equation}\label{eqn:lower_upper_bound}
\begin{aligned}
    \frac{F}{r_ip_i} &\leq \E[T(i,F)] \\ &= \frac{F}{r_ip_i} + \Delta \Char_{\{\alpha_i>0\}}  \frac{(1-\alpha_i) (1-p_i)}{p_i} \\ &\leq \frac{F}{r_ip_i} + \Delta \frac{1 - p_i}{p_i}.    
\end{aligned}
\end{equation}
To ensure $\E[T(i^*,F)]  \leq \E[T(j,F)]$ for all $j \in \cN \setminus \{ i^*\}$, it suffices to consider the upper bound of $E[T(i^*,F)]$ to be always smaller than the lower bound of $E[T(j,F)]$ for all $j \in \cN \setminus \{ i^*\}$ from \eqref{eqn:lower_upper_bound}, that is 
\begin{equation}\label{eqn:sufficient_inequality}
\begin{aligned}
\frac{F}{r_{i^*}p_{i^*}} + \Delta \frac{1-p_{i^*}}{p_{i^*}} \leq \min_{j \in \cN \setminus \{ i^*\}}\frac{F}{r_jp_j}.
\end{aligned}
\end{equation}
By definition of channel $h = \argmax_{j \in \cN \setminus \{i^*\}} r_jp_j$ we have $F/r_hp_h = \min_{j \in \cN \!\setminus\! \{ i^*\}}F/r_jp_j$. Then rearranging the second inequality in \eqref{eqn:sufficient_inequality} yields $F \geq H$ in (a).

When $F=k \Delta r_{i^*}$, we have $\alpha_{i^*} = 0$. Then from \eqref{eqn:time_static_policy}, the expected transfer time is simply $\E[T(i^*,F)] = \Delta\left(k/p_{i^*}\right) = F/r_{i^*}p_{i^*}$. Since $r_{i^*}p_{i^*} \geq r_jp_j$ for any $j\in\cN$, we have
$$\E[T(i^*,F)] = \frac{F}{r_{i^*}p_{i^*}} \leq \frac{F}{r_jp_j} \leq \E[T(j,F)]$$
for any $j \in \cN$, where the second inequality is from \eqref{eqn:lower_upper_bound}. Hence $i^*$ is the static optimal channel, that is, $i^* = i_{so}(F)$. This establishes (b), completing the proof.

\section{Proof of Corollary \ref{corollary:3.3}}\label{app:proof_coro_3.3}
For the file of size $F \in (k\Delta r_{i^*},(k+1)\Delta r_{i^*})$ with $k \in \N$, we have $r_{i^*}p_{i^*} > r_ip_i$ and $m_i \triangleq p_i/p_{i^*}$ for $i \in \cN \setminus \{i^*\}$. From \eqref{eqn:lower_upper_bound} in the proof of Proposition \ref{prop:threshold_optimal_policy}, we have $\E[T(i,F)] \leq F/r_ip_i+\Delta(1-p_i)/p_i$.
Moreover, from \eqref{eqn:time_static_policy} we have
\begin{equation}\label{eqn:lower_max-throughput_policy}
    \E[T(i^*,F)] > \Delta (k+1)(1/p_{i^*}-1) =  \Delta (k+1)(m_i/p_{i}-1).
\end{equation}
Therefore, the upper bound of the time ratio between channel $i$ and channel $i^*$ is shown as follows:
\begin{equation}\label{eqn:lower_bound_static_example}
\begin{split}
   \frac{\E[T(i,F)]}{\E[T(i^*,F)]} &< \frac{F/\Delta r_ip_i + (1-p_i)/p_i}{(k+1)(m_i/p_i-1)}.
\end{split}
\end{equation}
By definition of the static optimal channel and $\E[T(i_{so},F)]\!\leq\! \E[T(i^*,F)]$, we can get the result \eqref{eqn:lower_bound_static_general} by lower bounding \eqref{eqn:lower_bound_static_example} for channel $i \in \cN \setminus \{i^*\}$.

\section{Proof of Proposition \ref{prop:expected_time_dynamic_policy}}\label{app:proof_prop_4.1}
With our notation $E[T(\pi,s)]$ in mind. the Bellman equation for any fixed policy $\pi \in \Pi(F)$ (Proposition 4.2.3 in \cite{bertsekas2019reinforcement}) is shown as
\begin{equation}\label{eqn:Bellman_policy}
\begin{split}
     E[T(\pi,s)] =& \sum_{s' \in \cS} P_{s,s'}\left(\pi(s)\right)\left[c\left(s,\pi(s),s'\right) + E[T(\pi,s')]\right].
\end{split}
\end{equation}
The transition probability and cost function in section \ref{MDP setting} are defined as
\begin{equation*}
    P_{s,s}(i) = 1-p_i, ~~~~ P_{s,(s-\Delta r_i)^+}(i) = p_i, \forall s \in \mathcal{S} \!\setminus\! \{0\}, i \in \cN,
\end{equation*}
\begin{equation*}
    c\left(s,i,(s-\Delta r_i)^+\right) = \begin{cases} \min \{\Delta, s/r_i\} & s \in \mathcal{S} \!\setminus\! \{0\}, i \in \cN \\ 0 & s = 0, i \in \cN.  \end{cases}
\end{equation*}
Then by substituting $P_{s,s'}(i)$ and $c(s,i,s')$ with our transition probability and cost function defined above, \eqref{eqn:Bellman_policy} can be written as
\begin{equation*}
\begin{split}
     \E[T(\pi,F)] =& (1-p_{\pi(s)}) \left(\Delta + \E[T(\pi,F)]\right) \\ 
     &+ p_{\pi(s)} \left(\min\left\{\Delta ,\frac{s}{r_{\pi(s)}}\right\} + \E[T(\pi,F_2)]\right).
\end{split}
\end{equation*}

Recall that $F_n$ is the remaining file size right before the $n$-th successful file transmission given a policy $\pi$ and $F_{n+1} =F_n - \Delta r_{\pi(F_n)}$ for $n = 1,2,\cdots,|\pi|-1$, we can generalize this recursion to two adjacent states $F_n,F_{n+1} \in \cS$ in policy $\pi$ such that
\begin{equation}\label{eqn:recursion}
\begin{split}
  \E[T(\pi,F_n)] = &\Delta\frac{1-p_{\pi(F_n)}}{p_{\pi(F_n)}} + \min\left\{\Delta ,\frac{F_n}{r_{\pi(F_n)}}\right\} \\ &+   \E[T(\pi,F_{n+1})].    
\end{split}
\end{equation}
When $n = 1,2,\cdots,|\pi|-1$, the player fully spends $\Delta$ time in each successful transmission and the file transfer task has not been done yet $(F_n> \Delta r_{n})$, so that $\min\left\{\Delta ,F_n/r_{n}\right\} = \Delta$. For the last successful transmission $n=|\pi|$, we have 
$$\E[T(\pi,F_{|\pi|})] = \Delta\frac{1-p_{\pi_{|\pi|}}}{p_{\pi_{|\pi|}}} + \frac{F_{|\pi|}}{r_{\pi_{|\pi|}}}.$$ 
Thereby the expected transfer time with file size $F$ and policy $\pi$ can be recursively solved and the result in Proposition \ref{prop:expected_time_dynamic_policy} is derived.

\section{Proof of Proposition \ref{prop:lower_bound_coincide}}\label{app:proof_lemma_4.2}
From \eqref{eqn:recursion_closed_form} we have
\begin{equation*}
\begin{split}
   &\E[T(\pi,F)] \\ =& \sum_{n=1}^{|\pi|\!-\!1} \frac{\Delta r_{\pi(F_n)}}{r_{\pi(F_n)}p_{\pi(F_n)}} \!+\! \frac{\Delta\!\left(1\!-\!p_{\pi\left(F_{|\pi|}\right)}\right)\!r_{\pi\left(F_{|\pi|}\right)} \!+\! F_{|\pi|}p_{\pi\left(F_{|\pi|}\right)}}{r_{\pi\left(F_{|\pi|}\right)}p_{\pi\left(F_{|\pi|}\right)}} \\
   \geq& \sum_{n=1}^{|\pi|-1} \frac{\Delta r_{\pi(F_n)}}{r_{i^*}p_{i^*}} +\frac{\Delta\left(1\!-\!p_{\pi\left(F_{|\pi|}\right)}\right)r_{\pi\left(F_{|\pi|}\right)} \!+\! F_{|\pi|}p_{\pi\left(F_{|\pi|}\right)}}{r_{i^*}p_{i^*}} \\
   \geq& \frac{1}{r_{i^*}p_{i^*}}\left(\sum_{n=1}^{|\pi|-1}\Delta r_{\pi(F_n)} + F_{|\pi|}\right) = \frac{F}{r_{i^*}p_{i^*}},
\end{split}
\end{equation*}
where the first inequality comes from the fact that $r_{i^*}p_{i^*} \geq r_i p_i$ for all $i \in \cN$. The second inequality is from our definition of $F_{|\pi|}$ which implies that $F_{|\pi|} \leq \Delta r_{\pi_{|\pi|}}$.

When file size $F$ is an integer multiple of $\Delta r_{i^*}$, i.e., $F=k \Delta r_{i^*}$ for some $k \in \Z_+$, we have $\E[T(i^*,F)] = F/r_{i^*}p_{i^*}\leq \E[T(\pi^*,F)] \leq E[T(\pi,F)]$. Since $\E[T(i^*,F)] \geq \E[T(\pi^*,F)]$, we have $\E[T(\pi^*,F)] = F/r_{i^*}p_{i^*}$ as well, and the max-throughput policy coincides with the dynamic optimal policy.

\section{Proof of Corollary \ref{corollary:dynamic_lower_upper_bound}}\label{app:proof_coro_4.4}
We first define a suboptimal policy $\pi_{heu}$ and then use $E[T(\pi^*,F)] \leq E[T(\pi_{heu},F)]$ for our proof.
The suboptimal policy $\pi_{heu}$ is defined as sensing and accessing the max-throughput channel $i^*$ to transmit the file of size $F_1 = n\Delta r_{i^*}$, then following a static optimal policy for the remaining file of size $F_2 = F - n\Delta r_{i^*}$. $n$ is an integer chosen from $0$ to $k=\floor{F/\Delta r_{i^*}}$.
Then, the expected transfer time of the dynamic optimal policy is always smaller than that of the dynamic suboptimal policy, that is, 
\begin{equation*}
    \begin{split}
        \E[T(\pi^*,F)] &\leq \E[T(\pi_{heu},F)]\\
        &=\E[T(i^*,n\Delta r_{i^*})]+ \E[T(i_{so},F-n\Delta r_{i^*})]\\
        &\leq \Delta n/p_{i^*} +\min_{\substack{i \in \cN \!\setminus\!\{i^*\}}} \left\{ \frac{F-n\Delta r_{i^*}}{r_ip_i} +  \Delta\frac{1-p_i}{p_i}\right\} \\
        &= \!\!\min_{\substack{i \in \cN \!\setminus\!\{i^*\} }}\!\!\left\{\!\frac{F}{r_ip_i} + \Delta\left( \frac{1-p_i}{p_i} \!-\! n\frac{r_{i^*}p_{i^*}\!-\!r_ip_i}{r_ip_ip_{i^*}}\right)\!\right\}\!,
    \end{split}
\end{equation*}
where the second inequality comes from \eqref{eqn:time_static_policy}, and \eqref{eqn:lower_upper_bound} in the proof of Proposition \ref{prop:threshold_optimal_policy}. It shows monotonically decreasing in $n$ such that we can choose $n=k$ to get the smallest upper bound for $\E[T(\pi^*,F)]$. Moreover, we have $E[T(i^*,F)] > \Delta (k+1)(m_i/p_{i}-1)$ from \eqref{eqn:lower_max-throughput_policy}. Hence, the upper bound of the ratio $\E[T(\pi^*,F)]/\E[T(i^*,F)]$ in Corollary \ref{corollary:dynamic_lower_upper_bound} is proved.

For the lower bound of the ratio, by using Proposition \ref{prop:lower_bound_coincide} we have $E[T(\pi^*,F)] \geq F/r_{i^*}p_{i^*}$. Together with \eqref{eqn:static_lower_bound}, we have
\begin{equation*}
\begin{split}
     \frac{\E[T(\pi^*,F)]}{\E[T(i^*,F)]} &\geq \frac{F/r_{i^*}p_{i^*}}{F/r_{i^*}p_{i^*} + \Delta \Char_{\{\alpha_{i^*}>0\}}  (1\!-\!\alpha_{i^*}) (1\!-\!p_{i^*})/p_{i^*}} \\
     &\geq \frac{1}{1 + \Delta \Char_{\{\alpha_{i^*}>0\}}(1-p_{i^*})r_{i^*}/F},
\end{split}
\end{equation*}
where the second inequality comes from $1-\alpha_{i^*}\leq 1$. This completes the proof of Corollary \ref{corollary:dynamic_lower_upper_bound}.

\section{Proof of Theorem \ref{thm:regret_bound}}\label{app:proof_thm_5.1}
The proof is nearly the same as the analysis of Theorem 5.4 in Appendix G.B \cite{talebi2017stochastic}. Here we only give the main modifications for our setting.

The first modification comes from the definition of `arm'. In \cite{talebi2017stochastic}, each edge in the graph is treated as a different arm, that is, the status of each edge is observed and estimated separated. In our setting, each edge in the shortest path problem (Figure \ref{fig:equivalent_shortest_path}) represents one of $N$ channels such that each policy (path) may observe one channel multiple times. Then, some summation terms in the proof, previously were over all edges (e.g., (12), (13) in \cite{talebi2017stochastic}), are now over all $N$ channels. 

Second, the \textit{KL-SR} algorithm in \cite{talebi2017stochastic} for dynamic optimal policy works for a fixed source node, which can be interpreted as a fixed file size $F$. Our algorithm deals with the varying file size. With the boundness of file size $F_{\max}$, we only need to change parameter $H$ to be the longest policy length for maximum file size $F_{\max}$ (instead of fixed file size $F$), $\Delta_{\min}$ to be the smallest non-zero difference of expected transfer time between any sub-optimal policy and targeted optimal policy for file size in $(0,F_{\max}]$ (instead of fixed file size $F$) and $D$ to be the longest expected transfer time for file size $F_{\max}$ (instead of fixed file size $F$). Then, the proof will be carried over.

\end{document}